\documentclass[11pt]{article}
\usepackage{amssymb}
\usepackage{latexsym}
\usepackage{amsmath}
\usepackage{amsthm}
\newtheorem{theorem}{Theorem}
\newtheorem{lemma}[theorem]{Lemma}
\newtheorem{proposition}[theorem]{Proposition}
\newtheorem{corollary}[theorem]{Corollary}
\newtheorem{definition}[theorem]{Definition}

\newtheorem{observation}[theorem]{Observation}
\newtheorem{question}[theorem]{Question}
\newtheorem{remark}[theorem]{Remark}
\newcommand{\leg}[2]{\left(\begin{array}{c} #1 \\ \hline #2 \end{array}\right)}

\begin{document}
\date{Accepted for publication in Annals of Mathematics}
\title{Defining ${\mathbb Z}$ in ${\mathbb Q}$}

\author{Jochen Koenigsmann}
\maketitle

\begin{abstract}
We show that ${\mathbb Z}$ is definable in ${\mathbb Q}$
by a universal first-order formula in the language of rings.
We also present an $\forall\exists$-formula for ${\mathbb Z}$ in ${\mathbb Q}$
with just one universal quantifier.
We exhibit new diophantine subsets of ${\mathbb Q}$ like
the complement of the image of the norm map under a quadratic extension,
and we give an elementary proof for the fact
that the set of non-squares is diophantine.
Finally, we show that there is no existential formula for ${\mathbb Z}$ in ${\mathbb Q}$, provided one assumes a strong variant of the Bombieri-Lang Conjecture for varieties over ${\mathbb Q}$ with many ${\mathbb Q}$-rational points.\footnote{\noindent 2000
{\em Mathematics Subject Classification.}
Primary 11U05; Secondary 11R52, 11G35, 11U09. \newline
{\em Key words and phrases.}
Hilbert's Tenth Problem, diophantine set, undecidability, definability,
quaternion algebra, Bombieri-Lang Conjecture. \newline
The research on this paper started while the author enjoyed the hospitality
of the Max-Planck-Institut Bonn.}
\end{abstract}
\section{${\mathbb Z}$ is universally definable in ${\mathbb Q}$}
Hilbert's 10th problem was to find a general algorithm for deciding, given any $n$ and any polynomial $f\in {\mathbb Z} [x_1,\ldots ,x_n]$, whether or not $f$ has a zero in ${\mathbb Z}^n$.
Building on earlier work by Martin Davis, Hilary Putnam and Julia Robinson, Yuri Matiyasevich proved in 1970 that there can be no such algorithm. In particular, the existential first-order theory Th$_\exists ({\mathbb Z})$ of ${\mathbb Z}$ (in the language of rings ${\cal L}:=\{+, \cdot; 0,1\}$) is undecidable.
Hilbert's 10th problem over ${\mathbb Q}$, i.e., the question whether Th$_\exists ({\mathbb Q})$ is decidable, is still open.

If one had an {\bf existential} (or {\bf diophantine}) definition of ${\mathbb Z}$ in ${\mathbb Q}$ (i.e., a definition by an existential 1st-order ${\cal L}$-formula), then  Th$_\exists ({\mathbb Z})$ would be interpretable in Th$_\exists ({\mathbb Q})$, and the answer would, by Matiyasevich's Theorem, again be no. But it is still open whether ${\mathbb Z}$ is existentially definable in ${\mathbb Q}$, and, in fact, towards the end of the paper we provide a reason why it should not (Corollary \ref{notex}).

The earliest 1st-order definition of ${\mathbb Z}$ in ${\mathbb Q}$, due to Julia Robinson ([R]),
can be expressed by an $\forall\exists\forall$-formula of the shape
$$\phi (t):\;  \forall x_1\forall x_2\exists y_1\ldots \exists y_7\forall z_1\ldots\forall z_6\; f(t;x_1,x_2;y_1,\ldots ,y_7; z_1,\ldots,z_6 )=0$$
for some $f\in {\mathbb Z}[t;x_1,x_2;y_1,\ldots ,y_7; z_1,\ldots,z_6]$,
i.e., for any $t\in {\mathbb Q}$,
$$t\in{\mathbb Z}\mbox{ if and only if } \phi (t) \mbox{ holds in }{\mathbb Q}.$$
Recently, Bjorn Poonen ([P1]) managed to find an $\forall\exists$-definition with 2 universal and 7 existential quantifiers.
In this paper we present an $\forall$-definition of ${\mathbb Z}$ in ${\mathbb Q}$.
To search for such a creature is motivated by the following 
\medskip\\
{\bf Observation 0.} {\em If there is an existential definition of ${\mathbb Z}$ in ${\mathbb Q}$ then there is also a universal one.}
\medskip\\
{\em Proof:}
If ${\mathbb Z}$ is diophantine in ${\mathbb Q}$ then so is
$${\mathbb Q}\setminus {\mathbb Z} = \{ x\in {\mathbb Q}\mid\exists m,n,a,b\in{\mathbb Z}\mbox{ with } n\neq 0, \pm 1,\, am+bn=1\mbox{ and }m=xn\}.$$\qed
\begin{theorem}\footnote{In the meantime, Theorem 1 has been generalized
to arbitrary number fields $K$: the ring of integers of $K$
is universally definable in $K$ ([Pa]).}
There is a positive integer $n$ and a polynomial
$g\in {\mathbb Z}[t;x_1,\ldots ,x_n]$ such that, for any $t\in {\mathbb Q}$,
$$t\in {\mathbb Z}\mbox{ if and only if }
\forall x_1\ldots\forall x_n\in {\mathbb Q}\;\; g(t;x_1,\ldots ,x_n)\neq 0.$$
\end{theorem}
If one measures logical complexity in terms of the number of changes of quantifiers
then this is a definition of ${\mathbb Z}$ in ${\mathbb Q}$
of least possible complexity:
there is no quantifier-free definition of ${\mathbb Z}$ in ${\mathbb Q}$.
\begin{corollary}
${\mathbb Q}\setminus {\mathbb Z}$ is diophantine in ${\mathbb Q}$.
\end{corollary}
In more geometric terms, this says
\medskip\\
{\bf Corollary 2'.}
{\em There is a (not necessarily irreducible) affine variety $V$ over ${\mathbb Q}$ and a ${\mathbb Q}$-morphism $\pi:\, V\to {\mathbb A}^1$ such that the image of $V({\mathbb Q})$ is ${\mathbb Q}\setminus {\mathbb Z}$.}
\medskip\\
Together with the undecidability of Th$_\exists({\mathbb Z})$, Theorem 1 immediately implies
\begin{corollary}
{\em Th}$_{\forall\exists} ({\mathbb Q})$ is undecidable.
\end{corollary}
Here {\em Th}$_{\forall\exists} ({\mathbb Q})$
is the set of all sentences of the shape
$$\forall x_1\ldots\forall x_k\exists y_1\ldots\exists y_l\, \phi(x_1,\ldots ,x_k,y_1,\ldots,y_l),$$
where $\phi$ is a quantifier-free formula in the language of rings
${\cal L} =\{+,\cdot;0,1\}$, that is,
a boolean combination of polynomial equations and inequalities
between polynomials in $\mathbb{Z}[x_1,\ldots,x_k,y_1,\ldots,y_l]$.
Corollary 3 was proved conditionally, using a conjecture on elliptic curves, in [CZ].
Again, we can phrase this in more geometric terms:
\medskip\\
{\bf Corollary 3'.}
{\em There is no algorithm that decides on input a ${\mathbb Q}$-morphism $\pi:\, V\to W$ between affine ${\mathbb Q}$-varieties $V,W$
whether or not 
$\pi:\, V({\mathbb Q})\to W({\mathbb Q})$ is surjective.}
\medskip\\
{\bf Acknowledgement:}
Among many others, I would, in particular, like to thank Boris Zilber, Jonathan Pila, Marc Hindry and Joseph Silverman for very helpful discussions.
I am also most grateful to the anonymous referee
whose numerous suggestions have greatly improved the paper.
The now much shorter proofs of Proposition \ref{RnotQ}(a) and Corollary \ref{notex} are the referee's.
\section{The proof of Theorem 1}
Like all previous definitions of ${\mathbb Z}$ in ${\mathbb Q}$,
we use elementary facts on quadratic forms over $\mathbb{R}$ and $\mathbb{Q}_p$,
together with Hasse's Local-Global-Principle for quadratic forms.
What is new in our approach is the use of the Quadratic Reciprocity Law (e.g., in Proposition \ref{R_p} or \ref{RnotQ}) and, inspired by the model theory of local fields, the transformation of some existential formulas into universal formulas (Step 4).
A technical key trick is the existential definition of the Jacobson radical of certain rings (Step 3) which makes implicit use of so-called `rigid elements' as they occur, e.g., in [K].
\subsection*{{\em Step 1:} Diophantine definition of quaternionic semi-local rings \`{a} la Poonen}
The first step modifies Poonen's proof ([P1]), thus arriving at a formula for ${\mathbb Z}$ in ${\mathbb Q}$ which, like the formula in his Theorem 4.1, has 2 $\forall$'s followed by 7 $\exists$'s, but we managed to bring down the degree of the polynomial involved from 9244 to 8.
\begin{definition}
\label{Hab}
Let $\mathbb{P}$ be the set of rational primes and let $\mathbb{Q}_\infty :=\mathbb{R}$.

For $a,b\in {\mathbb Q}^\times$, let
\begin{itemize}
\item
$H_{a,b}:=  {\mathbb Q}\cdot 1\oplus {\mathbb Q}\cdot \alpha\oplus {\mathbb Q}\cdot \beta\oplus {\mathbb Q}\cdot \alpha\beta$ be
the quaternion algebra over ${\mathbb Q}$ with multiplication defined by
$\alpha^2 =a$, $\beta^2 =b$ and $\alpha\beta = - \beta\alpha$,
\item
$\Delta_{a,b}:=\{p\in {\mathbb P}\cup\{\infty\}\mid H_{a,b}\otimes {\mathbb Q}_p\not\cong M_2({\mathbb Q}_p)\}$
the set of primes (including $\infty$) where $H_{a,b}$ does not split locally
--- $\Delta_{a,b}$ is always finite, and $\Delta_{a,b}=\emptyset$ iff $a\in N(b)$, i.e., $a$ is in the image of the norm map ${\mathbb Q} (\sqrt{b})\to {\mathbb Q}$,
\item
$S_{a,b}:= \{2x_1\in{\mathbb Q}\mid \exists x_2,x_3,x_4\in{\mathbb Q}:\,x_1^2-ax_2^2-bx_3^2+abx_4^2=1\}$ the set of traces of norm-$1$ elements of $H_{a,b}$, and
\item
$T_{a,b}:= S_{a,b}+S_{a,b}$ -- note that $T_{a,b}$ is an existentially defined subset of ${\mathbb Q}$. Here we deviate from Poonen's terminology: his $T_{a,b}$ is $S_{a,b}+ S_{a,b}+\{0,1,\ldots ,2309\}$.
\end{itemize}
\end{definition}
For each $p\in\mathbb{P}\cup\{\infty\}$,
we can similarly define $S_{a,b}(\mathbb{Q}_p)$ and $T_{a,b}(\mathbb{Q}_p)$
by replacing $\mathbb{Q}$ by $\mathbb{Q}_p$.

For each $p\in\mathbb{P}$,
we will denote the $p$-adic valuation on $\mathbb{Q}$ or on $\mathbb{Q}_p$ by $v_p$,
and the assoicated residue map by $\phi_p:\,\mathbb{Z}_{(p)}\to\mathbb{F}_p$
resp. $\phi_p:\,\mathbb{Z}_p\to\mathbb{F}_p$.

An explicit criterion for checking whether or not an element $p\in {\mathbb P}\cup\{\infty\}$ belongs to $\Delta_{a,b}$, is given in the following
\begin{observation}
\label{pinDelta}
Assume $a,b\in {\mathbb Q}^\times$ and $p\in {\mathbb P}\cup\{\infty\}$.
Then $p\in\Delta_{a,b}$ if and only if
\begin{description}
\item[for ${\bf p=2}$:]
After multiplying by suitable rational squares and integers $\equiv 1\!\!\!\!\mod 8$ and, possibly, swapping $a$ and $b$, the pair $(a,b)$ is one of the following:
$$\begin{array}{lllll}
(2,3) & (3,3) & (5,6) & (6,6) & (15,15)\\
(2,5) & (3,10) & (5,10) & (6,15) & (15,30)\\
(2,6) & (3,15) & (5,30) & (10,30) & (30,30)\\
(2,10) & & & &
\end{array}$$
\item[for ${\bf 2\neq p\in {\mathbb P}}$:]
$$\begin{array}{l}
v_p(a)\mbox{ is odd, }v_p(b)\mbox{ is even, and }\left(\frac{bp^{-v_p(b)}}{p}\right) =-1\mbox{, or}\\
v_p(a)\mbox{ is even, }v_p(b)\mbox{ is odd, and } \left(\frac{ap^{-v_p(a)}}{p}\right) =-1\mbox{ or}\\
v_p(a)\mbox{ is odd, }v_p(b)\mbox{ is odd, and } \left(\frac{-abp^{-v_p(ab)}}{p}\right)=-1.\\
\end{array}$$
\item[for ${\bf p=\infty}$:]
$a<0$ and $b<0$.
\end{description}
\end{observation}
{\em Proof:}
This is an immediate translation of the computation of the Hilbert symbol
$(a,b)_p$ (which is $1$ or $-1$ depending on whether or not $p\in\Delta_{a,b}$)
as in Theorem 1 of Ch.III in [Se]:\\
For finite odd $p$ and $a=p^\alpha u$ and $b=p^\beta v$ (with $u,v$ $p$-adic units) the formula is
$$(a,b)_p= (-1)^{\alpha\beta\epsilon (p)}\leg{u}{p}^\beta\leg{v}{p}^\alpha,$$
where $\epsilon(p):=\frac{p-1}{2}\!\!\mod 2$.\\
For $p=2$, the formula is
$$(a,b)_2=(-1)^{\epsilon (u)\epsilon (v) +\alpha\omega (v) +\beta\omega (u)},$$
where $\omega (u):=\frac{u^2-1}{8}\!\!\mod 2$.\\
For $p=\infty$ the statement is obvious.
\qed
\begin{proposition}
\label{Tab}
For any $a,b\in {\mathbb Q}^\times$,
$$T_{a,b}=\bigcap_{p\in\Delta_{a,b}}{\mathbb Z}_{(p)},$$
where ${\mathbb Z}_{(\infty)}:=\{ x\in{\mathbb Q}\mid -4\leq x\leq 4\}$. 
\end{proposition}
Here and throughout the rest of the paper, we use the following
\medskip\\
{\bf Convention: }{\em Given an empty collection of subsets of $\mathbb{Q}$,
the intersection is $\mathbb{Q}$.}
\medskip

{\em Proof:}
For each $p\in {\mathbb P}$, let $$U_p:=\{ s\in {\mathbb F}_p\mid x^2-sx +1\mbox{ is irreducible over }{\mathbb F}_p\}$$.
We shall use the following
\medskip\\
{\bf Facts} {\em For any $a,b\in{\mathbb Q}^\times$ and for any $p\in {\mathbb P}$:
\begin{description}
\item[(a)]
If $p\not\in\Delta_{a,b}$ then $S_{a,b}({\mathbb Q}_p) ={\mathbb Q}_p$.
\item[(b)]
If $p\in\Delta_{a,b}$ then $\phi_p^{-1}(U_p)\subseteq S_{a,b} ({\mathbb Q}_p)\subseteq {\mathbb Z}_p$.
\item[(c)]
$S_{a,b}({\mathbb R})=\left\{ \begin{array}{ll}
{\mathbb R} & \mbox{for }a>0\mbox{ or }b>0\\
 \mbox{{\em [}} -2,2 \mbox{{\em ]}} & \mbox{for }a,b<0.
\end{array}\right.$
\item[(d)]
If $p>11$ then ${\mathbb F}_p = U_p + U_p$.
\item[(e)]
$S_{a,b}({\mathbb Q}) ={\mathbb Q}\cap\bigcap_{p\in\Delta_{a,b}}S_{a,b}({\mathbb Q}_p)$.
\end{description}}
{\bf (a)} and {\bf (b)} are [P1], Lemma 2.1, {\bf (c)} is a straightforward computation, {\bf (d)} is [P1], Lemma 2.3, and {\bf (e)} is a special case of the Hasse-Minkowski local-global principle for representing rationals by quadratic forms.

{\bf (b)} and {\bf (c)} immediately give the inclusion $T_{a,b}\subseteq \bigcap_{p\in\Delta_{a,b}}{\mathbb Z}_{(p)}$.

To prove the converse inclusion $T_{a,b}\supseteq \bigcap_{p\in\Delta_{a,b}}{\mathbb Z}_{(p)}$,
let us first compute $U_p$ for the primes $p\leq 11$:
$$\begin{array}{lll}
U_2 & = & \{1\}\\
U_3 & = & \{0\}\\
U_5 & = & \{1,4\}\\
U_7 & = & \{0,3,4\}\\
U_{11} & = & \{0,1,5,6,10\}.
\end{array}$$

For each $p\in {\mathbb P}\cup\{\infty\}$ define $V_p\subseteq {\mathbb Z}_p$ as follows:
$$V_p = \left\{\begin{array}{ll}
\phi_2^{-1}(U_2)\cup (4+8{\mathbb Z}_2) & \mbox{for }p=2\\
\phi_p^{-1}(U_p)\cup [(\pm 2 + p{\mathbb Z}_p)\setminus (\pm 2 + p^2{\mathbb Z}_p)] & \mbox{for }3\leq p\leq 11\\
\phi_p^{-1}(U_p) & \mbox{for }11< p\in{\mathbb P}\\
 \mbox{[} -2,2 \mbox{]} & \mbox{for } p=\infty.
\end{array}\right.$$
(We define ${\mathbb Z}_\infty$ to be the real interval $[-4,4]\subseteq {\mathbb R}$.)\\
By Fact {\bf (b)}, Fact {\bf (c)}, Observation \ref{pinDelta} together with an easy direct calculation in the cases $p=3,5,7,11$ and, for $p=2$, by the table below, one always has
$$V_p\subseteq S_{a,b}({\mathbb Q}_p)\mbox{ and, for }p\neq\infty,\,V_p\mbox{ is open}.$$

The table for $p=2$ lists those pairs $(a,b)$ with $(a,b)_2=-1$
as in Observation \ref{pinDelta}, and gives, in each case,
$$4+8{\mathbb Z}_2\subseteq S_{a,b}({\mathbb Q}_2)$$
by assuming that we are given
$x_1\in 2+8\mathbb{Z}_2$ or $x_1\in 6+8\mathbb{Z}_2$
(which is equivalent to $2x_1\in 4+8\mathbb{Z}_2$)
and by specifying elements $x_2, x_3$ and $x_4$ which guarantee that
$$-ax_2^2-bx_3^2+abx_4^2\equiv_2 1-x_1^2\equiv_2 -3\!\!\mod 8\mathbb{Z}_2.$$
Multiplying $x_2^2, x_3^2, x_4^2$ by a suitable common element from $1+8{\mathbb Z}_2\subseteq ({\mathbb Q}_2^\times)^2$, makes then sure that $2x_1\in S_{a,b}({\mathbb Q}_2)$.
$$\begin{array}{ccccc}
(a,b) & \vline\,\vline & x_2 & x_3 & x_4\\ \hline\hline
(2,3) & \vline\,\vline & 0 & 1 & 0\\ 
(2,5) & \vline\,\vline & 2 & 1 & 1\\
(2,6) & \vline\,\vline & 0 & 1 & \frac{1}{2}\\
(2,10) & \vline\,\vline & 2 & 0 & \frac{1}{2}\\ \hline
(3,3) & \vline\,\vline & 1 & 0 & 0\\
(3,10) & \vline\,\vline & 1 & 0 & 0\\
(3,15) & \vline\,\vline & 1 & 0 & 0\\ \hline
(5,6) & \vline\,\vline & 1 & 1 & 0\\
(5,10) & \vline\,\vline & 1 & 0 & 1\\
(5,30) & \vline\,\vline & 1 & 1 & 0\\ \hline
(6,6) & \vline\,\vline & \frac{1}{2} & \frac{1}{2} & 0\\
(6,15) & \vline\,\vline & 1 & 1 & 0\\ \hline
(10,30) & \vline\,\vline & 0 & 1 & \frac{1}{10}\\ \hline
(15,15) & \vline\,\vline & 1 & 0 & \frac{2}{15}\\
(15,30) & \vline\,\vline & 1 & 1 & \frac{1}{15}\\ \hline
(30,30) & \vline\,\vline & 1 & 1 & \frac{1}{30}
\end{array}$$

Fact {\bf (d)} and another elementary case-by-case-check for $p\leq 11$ shows
that for any $p\in{\mathbb P}\cup \{\infty\}$
$${\mathbb Z}_p = V_p + V_p.$$
Now pick $t\in\bigcap_{p\in\Delta_{a,b}} {\mathbb Z}_{(p)}$.
For each $p\in\Delta_{a,b}$, there is some $s_p\in {\mathbb Z}_p$
such that $s_p, t-s_p\in V_p$.

If $t=\pm 4$ then, clearly, $t=\pm 2 \pm 2\in S_{a,b}+S_{a,b}=T_{a,b}$.

If $t\neq\pm 4$ and $\infty\in\Delta_{a,b}$
we can choose $s_\infty\in {\mathbb Z}_\infty = [-4,4]\subseteq {\mathbb R}$ such that
$s_\infty, t-s_\infty\in ]-2,2[$.
Now approximate the finitely many $s_p\in {\mathbb Z}_p$ ($p\in\Delta_{a,b}$)
by a single $s\in {\mathbb Q}$ such that
$$s-s_p\in\left\{\begin{array}{ll}
8 {\mathbb Z}_2 & \mbox{if }p=2\\
p^2{\mathbb Z}_p & \mbox{if }3\leq p\leq 11\\
p{\mathbb Z}_p & \mbox{if }11<p\in {\mathbb P}\\
 ]-\epsilon,\epsilon[ & \mbox{if }p=\infty
\end{array}\right.$$
where $\epsilon = \min \{\mid 2\pm s_\infty\mid,\,\mid 2\pm (t-s_\infty)\mid\}$.
This guarantees that for all $p\in\Delta_{a,b}$
$$s, t-s\in V_p\subseteq S_{a,b}({\mathbb Q}_p),$$
and hence, by Fact {\bf (e)}, that $s,t-s\in S_{a,b} = S_{a,b}({\mathbb Q})$.\qed
\medskip

One then obtains an $\forall\exists$-definition of ${\mathbb Z}$ in ${\mathbb Q}$ from the fact that
$${\mathbb Z}=\bigcap_{l\in{\mathbb P}}{\mathbb Z}_{(l)}=\bigcap_{a,b>0} T_{a,b}$$
as in [P1], Theorem 4.1.
With our simplified $T_{a,b}$, the formula now becomes, for any $t\in {\mathbb Q}$,
$$t\in {\mathbb Z} \Longleftrightarrow
\begin{array}{l}
\forall a,b\,\exists x_1, x_2, x_3, x_4, y_2, y_3, y_4\\
(a+x_1^2+x_2^2+x_3^2+x_4^2)\cdot (b+x_1^2+x_2^2+x_3^2+x_4^2)\cdot \\
\mbox{[} (x_1^2-ax_2^2-bx_3^2+abx_4^2-1)^2 +\\
+ \,((t-2x_1)^2-4ay_2^2-4by_3^2+4aby_4^2 -4)^2\mbox{]} \,=\,0
\end{array}$$
\subsection*{{\em Step 2:} Towards a uniform diophantine definition of all ${\mathbb Z}_{(p)}$'s in ${\mathbb Q}$}
We will present a diophantine definition for the local rings ${\mathbb Z}_{(p)}={\mathbb Z}_p\cap {\mathbb Q}$ depending on the congruence of the prime $p$ modulo $8$, and involving $p$ (and if $p\equiv 1\!\!\mod 8$ an auxiliary prime $q$) as a parameter. However, since in any first-order definition of a subset of ${\mathbb Q}$ we can only quantify over the elements of ${\mathbb Q}$, and not, e.g., over all primes, we will allow arbitrary non-zero rational numbers $p$ and $q$ as parameters in the following definition. 
\begin{definition}
\label{R}
For $p,q\in {\mathbb Q}^\times$, let
\begin{itemize}
\item 
$R_p^{[3]} :=T_{-1,-p}+T_{2,-p}$
\item 
$R_p^{[5]} :=T_{-2,-p}+T_{2,-p}$
\item 
$R_p^{[7]} :=T_{-1,-p}+T_{-2,p}$
\item
$R_{p,q}^{[1]} :=T_{-2p,q}+T_{2p,q}$
\end{itemize}
\end{definition}
\begin{remark}
\label{T+T}
\begin{enumerate}
\item[(a)]
For any $a,b,c,d\in\mathbb{Q}^\times$ with at least one of them positive,
$$T_{a,b}+T_{c,d}=
\bigcap_{l\in\Delta_{a,b}}\mathbb{Z}_{(l)}+\bigcap_{l\in\Delta_{c,d}}\mathbb{Z}_{(l)} 
= \bigcap_{l\in\Delta_{a,b}\cap\Delta_{c,d}}\mathbb{Z}_{(l)}$$
\item[(b)]
The $R$'s are are existentially defined, uniform in $p$ and $q$,
so that for $k=3,5$ or $7$ the sets
$$\{ (p,x)\in\mathbb{Q}^\times\times\mathbb{Q}\mid x\in R_p^{[k]}\}$$
and the set
$$\{ (p,q,x)\in\mathbb{Q}^\times\times\mathbb{Q}^\times\mathbb{Q}\mid x\in R_{p,q}^{[1]}\}$$
are diophantine.
\end{enumerate}
\end{remark}
{\em Proof:}
(a) The first equation is from Proposition \ref{Tab}.
For the second equation, the inclusion `$\subseteq$' is obvious.
For `$\supseteq$', assume $x\in\bigcap_{l\in\Delta_{a,b}\cap\Delta_{c,d}}\mathbb{Z}_{(l)}$.
By approximation, there is $y\in\mathbb{Q}$ such that
$$y\in\left\{\begin{array}{cl}
x+l\mathbb{Z}_{(l)} & \mbox{for }l\in\Delta_{c,d}\\
\mathbb{Z}_{(l)} & \mbox{for }l\in\Delta_{a,b}\setminus\Delta_{c,d}
\end{array}\right.$$
Then $y\in\bigcap_{l\in\Delta_{a,b}}\mathbb{Z}_{(l)}$
and $x-y\in\bigcap_{l\in\Delta_{c,d}}\mathbb{Z}_{(l)}$,
so that $x=y+(x-y)\in
\bigcap_{l\in\Delta_{a,b}}\mathbb{Z}_{(l)}+\bigcap_{l\in\Delta_{c,d}}\mathbb{Z}_{(l)}$.
\medskip\\
(b) This is immediate from definitions \ref{Hab} and \ref{R}.\qed
\begin{definition}
\label{P(p)}
\begin{enumerate}
\item[(a)]
For $k=1,3,5$ or $7$, define
${\mathbb P}^{[k]}:=\{ l\in {\mathbb P}\mid l\equiv k\!\!\mod 8\}$
\item[(b)] 
For $p\in {\mathbb Q}^\times$, define
\begin{itemize}
\item 
${\mathbb P}(p):=\{ l\in {\mathbb P}\mid v_l(p)\mbox{ is odd}\}$
\item 
${\mathbb P}^{[k]}(p):={\mathbb P}(p)\cap {\mathbb P}^{[k]}$, where $k=1,3,5$ or $7$.
\end{itemize}
\end{enumerate}
\end{definition}
\begin{proposition}
\label{R_p}
\begin{enumerate}
\item[(a)] 
${\mathbb Z}_{(2)} = T_{3,3}+T_{2,5}$
\item[(b)]
Suppose that $k=3,5$ or $7$. Then, for $p\in\mathbb{Q}^\times$, 
$$R_p^{[k]}=\left\{\begin{array}{ll}
\bigcap_{l\in {\mathbb P}^{[k]}(p)} {\mathbb Z}_{(l)} & \mbox{if } p\equiv k \;(\!\!\!\mod 8\mathbb{Z}_{(2)}) \\
\bigcap_{l\in {\mathbb P}^{[k]}(p)} {\mathbb Z}_{(l)}\mbox{ or }\bigcap_{l\in {\mathbb P}^{[k]}(p)\cup\{2\}} {\mathbb Z}_{(l)} & \mbox{otherwise}
\end{array}\right.$$
(As before, $\bigcap_{l\in\emptyset} {\mathbb Z}_{(l)} = {\mathbb Q}$.)

In particular, if $p$ is a prime and $p\equiv k \!\!\mod 8$ then ${\mathbb Z}_{(p)} = R_p^{[k]}$.
\item[(c)]
For $p,q\in {\mathbb Q}^\times$ with $p\equiv 1\;(\!\!\!\mod 8\mathbb{Z}_{(2)})$ and $q\equiv 3\;(\!\!\!\mod 8\mathbb{Z}_{(2)})$,
$$R_{p,q}^{[1]} = \bigcap_{l\in {\mathbb P}(p,q)} {\mathbb Z}_{(l)}$$
where ${\mathbb P}(p,q):= \Delta_{-2p,q}\cap\Delta_{2p,q}$.

In particular, if $p$ is a prime $\equiv 1\!\!\mod 8$ and $q$ is a prime $\equiv 3\!\!\mod 8$ with $\leg{p}{q} =-1$ then ${\mathbb Z}_{(p)}=R_{p,q}^{[1]}$.
\end{enumerate}
\end{proposition}
{\em Proof:}
(a) By Observation \ref{pinDelta}, $\Delta_{3,3} =\{ 2,3\}$ and $\Delta_{2,5}=\{ 2,5\}$, hence, by Remark \ref{T+T}(a),
$$ T_{3,3}+T_{2,5} =\bigcap_{l\in\Delta_{3,3}\cap\Delta_{2,5}}{\mathbb Z}_{(l)} = {\mathbb Z}_{(2)}.$$
(b) First assume $p\in {\mathbb Q}^\times$ with $p\equiv 3\;(\!\!\!\mod 8\mathbb{Z}_{(2)})$. Then, by Observation \ref{pinDelta},
$$\begin{array}{ccl}
\Delta_{-1,-p}\cap {\mathbb P} & = & {\mathbb P}^{[3]}(p)\cup {\mathbb P}^{[7]}(p)\\
\Delta_{2,-p} & = & {\mathbb P}^{[3]}(p)\cup {\mathbb P}^{[5]}(p)\cup \{ 2\},
\end{array}$$
so $\Delta_{-1,-p}\cap\Delta_{2,-p}= {\mathbb P}^{[3]}(p)$,
and, by Remark \ref{T+T}(a),
$$R_p^{[3]}:=T_{-1,-p}+T_{2,-p}=\bigcap_{l\in \Delta_{-1,-p}\cap\Delta_{2,-p}}{\mathbb Z}_{(l)} = \bigcap_{l\in {\mathbb P}^{[3]}(p)} {\mathbb Z}_{(l)}.$$
If $p\not\equiv 3\;(\!\!\!\mod 8\mathbb{Z}_{(2)})$, the only possible additional prime is $2$
(e.g. if $p\equiv 5\;(\!\!\!\mod 8\mathbb{Z}_{(2)})$).

If $p\equiv 5\;(\!\!\!\mod 8\mathbb{Z}_{(2)})$ then, again by Observation \ref{pinDelta},
$$\begin{array}{ccl}
\Delta_{-2,-p}\cap {\mathbb P} & = & {\mathbb P}^{[5]}(p)\cup {\mathbb P}^{[7]}(p)\\
\Delta_{2,-p} & = & {\mathbb P}^{[3]}(p)\cup {\mathbb P}^{[5]}(p)\cup\{ 2\},
\end{array}$$
so $\Delta_{-2,-p}\cap\Delta_{2,-p}= {\mathbb P}^{[5]}(p)$,
and $$R_p^{[5]}:=T_{-2p,-p}+T_{2p,-p}=\bigcap_{l\in \Delta_{-2,-p}\cap\Delta_{2,-p}}{\mathbb Z}_{(l)} =
\bigcap_{l\in {\mathbb P}^{[5]}(p)} {\mathbb Z}_{(l)}.$$
Again, the prime $2$ (and no other prime) may or may not enter if $p\not\equiv 5\;(\!\!\!\mod 8\mathbb{Z}_{(2)})$.

Finally, if $p\equiv 7\;(\!\!\!\mod 8\mathbb{Z}_{(2)})$ then, again by Observation \ref{pinDelta},
$$\begin{array}{ccl}
\Delta_{-1,-p}\cap {\mathbb P} & = & {\mathbb P}^{[3]}(p)\cup {\mathbb P}^{[7]}(p)\\
\Delta_{-2,p}\cap {\mathbb P} & = & {\mathbb P}^{[5]}(p)\cup {\mathbb P}^{[7]}(p)\cup\{ 2\},
\end{array}$$
so $\Delta_{-1,-p}\cap\Delta_{-2,p}= {\mathbb P}^{[7]}(p)$,
and $$R_p^{[7]}:=T_{-p,-p}+T_{2p,p}=\bigcap_{l\in \Delta_{-1,-p}\cap\Delta_{-2,p}}{\mathbb Z}_{(l)} =\bigcap_{l\in {\mathbb P}^{[7]}(p)} {\mathbb Z}_{(l)}.$$
As before, $2$ may enter if $p\not\equiv 7\;(\!\!\!\mod 8\mathbb{Z}_{(2)})$.

(c) The first statement is immediate from Remark \ref{T+T}(a).
For the `in particular', assume $p$ and $q$ are primes
with $p\equiv 1\!\!\mod 8$, $q\equiv 3\!\!\mod 8$
and $\left(\frac{p}{q}\right)=-1$.
Then, by quadratic reciprocity, $\left(\frac{q}{p}\right)=-1$,
and so, from Observation \ref{pinDelta},
$\Delta_{-2p,q}=\{p\}$ and $\Delta_{2p,q}=\{ p,q\}$.
Hence $R_{p,q}^{[1]}=\mathbb{Z}_{(p)}$.\qed
\begin{corollary}
\label{ZviaRs} 
$${\mathbb Z} ={\mathbb Z}_{(2)}\cap\bigcap_{p,q\in{\mathbb Q}^\times}(R_p^{[3]}\cap R_p^{[5]}\cap R_p^{[7]}\cap R_{p,q}^{[1]})$$
\end{corollary}
{\em Proof:}
By Remark \ref{T+T}(a), all $R$'s on the right hand side
are semilocal subrings of $\mathbb{Q}$ containing $\mathbb{Z}$.
On the other hand, by the `in particular' parts of the proposition,
for each prime $p$,  the right hand side is contained in $\mathbb{Z}_{(p)}$;
note that for $p\equiv 1\!\!\mod 8$ one always finds a prime $q\equiv 3\!\!\mod 8$
such that $q$ is congruent to a non-square $\!\!\mod p$.\qed
\subsection*{{\em Step 3:} An existential definition for the Jacobson radical}
We will show that, for some rings $R$ occuring in Proposition \ref{R_p}, the Jacobson radical $J(R)$ can be defined by an existential formula. This will also give rise to new diophantine predicates in ${\mathbb Q}$.
\begin{definition}
For $a,b,c\in\mathbb{Q}^\times$ we define
\begin{itemize}
\item 
$T_{a,b}^\times :=\{ u\in T_{a,b}\mid\exists v\in T_{a,b}\mbox{ with }uv=1\}$
\item 
$I_{a,b}^c:= c\cdot\mathbb{Q}^2\cdot T_{a,b}^\times\cap (1-\mathbb{Q}^2\cdot T_{a,b}^\times )$
\item 
$J_{a,b}:= (I_{a,b}^a + I_{a,b}^a)\cap (I_{a,b}^b + I_{a,b}^b)$
\end{itemize}
\end{definition}
Note that the set $\{ (a,b,x)\in \mathbb{Q}^\times\times\mathbb{Q}^\times\times\mathbb{Q}\mid x\in J_{a,b}\}$
is diophantine.
\begin{lemma}
\label{Jab}
Assume $a,b,c\in\mathbb{Q}^\times$. Then
\begin{enumerate}
\item[(a)]
$T_{a,b}^\times = \left\{\begin{array}{ll}
\bigcap_{l\in\Delta_{a,b}}\mathbb{Z}_{(l)}^\times & \mbox{if }\infty\not\in\Delta_{a,b}\\
([-4,-\frac{1}{4}]\cup[\frac{1}{4}, 4]) \cap\bigcap_{l\in\Delta_{a,b}\setminus\{\infty\}}{\mathbb Z}_{(l)}^\times
& \mbox{if }\infty\in\Delta_{a,b}\end{array}\right.$
\item[(b)]
$I_{a,b}^c = \{ 0\}\cup\left\{ y\in\mathbb{Q}^\times\vline
\begin{array}{l}v_l(y)\mbox{is odd and positive for all }l\in\Delta_{a,b}\cap\mathbb{P}(c)\mbox{ and}\\
v_l(y),v_l(1-y)\mbox{are even for all }l\in\Delta_{a,b}\setminus (\mathbb{P}(c)\cup\{\infty\})
\end{array}\right\}$\footnote{Here we adopt the convention that $\infty$ is even (to include the case that $y=1$
which can only happen when $\Delta_{a,b}\cap\mathbb{P}(c)=\emptyset$,
a case that will never be used later).}
\item[(c)]
$I_{a,b}^c+ I_{a,b}^c = \bigcap_{l\in\Delta_{a,b}\cap {\mathbb P}(c)}l\,{\mathbb Z}_{(l)},$
\item[(d)]
$J_{a,b} =\bigcap_{l\in\Delta}l\mathbb{Z}_{(l)}$,
where $\Delta =\left\{\begin{array}{ll}
\Delta_{a,b}\setminus\{ 2,\infty\} & \mbox{if }2\in\Delta_{a,b}\mbox{ and }
v_2(a), v_2(b)\mbox{ are even}\\
\Delta_{a,b}\setminus\{\infty\} & \mbox{else}
\end{array}\right.$ 
\end{enumerate}
In particular, if $\infty\not\in\Delta_{a,b}$
then $T_{a,b}^\times$ is the group of units of the ring $T_{a,b}$
and, if also $2\not\in\Delta_{a,b}$ or at least one of $v_2(a), v_2(b)$ is odd,
$J_{a,b}$ is the Jacobson radical of $T_{a,b}$.
\end{lemma} 
{\em Proof:}
(a) This is an immediate consequence of Proposition \ref{Tab}.

(b) `$\subseteq$':
By weak approximation,
$$\mathbb{Q}^2\cdot T_{a,b}^\times =\{ 0\}\cup\bigcap_{l\in\Delta_{a,b}\setminus\{\infty\}}v_l^{-1}(2\mathbb{Z}).$$
So if $y\in I_{a,b}^c\setminus\{ 0\}$ and $l\in\Delta_{a,b}\cap\mathbb{P}(c)$
then $v_l(y)$ is odd.
On the other hand, $1-y\in\mathbb{Q}^2\cdot T_{a,b}^\times$,
so $v_l(1-y)$ is even.
By the ultrametric inequality, this is only possible when $v_l(y)>0$.
If, on the other hand, $l\in\Delta_{a,b}\setminus (\mathbb{P}(c)\cup\{\infty\})$ then $v_l(y)$ and $v_l(1-y)$ are even.\\
`$\supseteq$':
Clearly, $0\in I_{a,b}^c$.
Now assume $y\in\mathbb{Q}^\times$
such that, for all $l\in\Delta_{a,b}\cap\mathbb{P}(c)$,
$v_l(y)$ is positive and odd.
Then $c^{-1}y\in \bigcap_{l\in\Delta_{a,b}\cap\mathbb{P}(c)}v_l^{-1}(2\mathbb{Z})$
and $1-y\in\bigcap_{l\in\Delta_{a,b}\cap\mathbb{P}(c)}\mathbb{Z}_{(l)}^\times\subseteq
\bigcap_{l\in\Delta_{a,b}\cap\mathbb{P}(c)}v_l^{-1}(2\mathbb{Z})$.

If we assume that $v_l(y)$ and $v_l(1-y)$ are even for all
$l\in\Delta^\prime:=\Delta_{a,b}\setminus (\mathbb{P}(c)\cup\{\infty\})$
then both $c^{-1}y$ and $1-y$ lie in $\bigcap_{l\in\Delta^\prime}v_l^{-1}(2\mathbb{Z})$.

So with both assumptions we see that both $c^{-1}y$ and $1-y$
lie in $$\bigcap_{l\in\Delta_{a,b}\setminus\{\infty\}}v_l^{-1}(2\mathbb{Z})
\subseteq \mathbb{Q}^2\cdot T_{a,b}^\times$$.

(c) For any prime $l$, any $x\in\mathbb{Q}$ with $v_l(x)>0$ can be written as the sum of two elements of odd positive value.
And any $x\in\mathbb{Q}$ can be written as the sum of two elements $y_1$ and $y_2$ such that $v_l(y_i)$ and $v_l(1-y_i)$ are both even for both $i=1,2$:
choose $y_1$ of even value $<\min \{ 0, v_l(x)\}$ and let $y_2=x-y_1$; then
$v_l(1-y_1)=v_l(y_1)=v_l(y_2)=v_l(1-y_2)$.
Hence the claim follows by approximation. 

(d) By definition, $J_{a,b} = (I_{a,b}^a + I_{a,b}^a)\cap (I_{a,b}^b+I_{a,b}^b)$,
so, from (c),
$$J_{a,b} =\bigcap_{l\in\Delta_{a,b}\cap\mathbb{P}(a)}l\mathbb{Z}_{(l)}
\cap \bigcap_{l\in\Delta_{a,b}\cap\mathbb{P}(b)}l\mathbb{Z}_{(l)}
= \bigcap_{l\in\Delta_{a,b}\cap (\mathbb{P}(a)\cup\mathbb{P}(b))}l\mathbb{Z}_{(l)},$$
where the second equality is, again, by weak approximation.
But now, from Observation \ref{pinDelta},
$$\Delta_{a,b}\cap (\mathbb{P}(a)\cup\mathbb{P}(b))
=\left\{\begin{array}{ll}
\Delta_{a,b}\setminus\{ 2,\infty\} & \mbox{if }2\in\Delta_{a,b}\mbox{ and }
v_2(a), v_2(b)\mbox{ are even}\\
\Delta_{a,b}\setminus\{\infty\} & \mbox{else}
\end{array}\right.$$
\qed
\medskip

Before we give the existential definition of the Jacobson radical $J(R)$ for some of the rings $R$ defined in Step 2
(Corollary \ref{Phi_k} and Proposition \ref{RnotQ} below)
we require another easy Lemma:
\begin{lemma}
\label{J+J}
Let $a,b,c,d\in\mathbb{Q}^\times$,
at least one of which positive
and at least one of which with odd dyadic value.
Let $\Delta:=\Delta_{a,b}\cap\Delta_{c,d}$ and let
$R=\bigcap_{l\in\Delta}\mathbb{Z}_{(l)}$.
Then $$J_{a,b} + J_{c,d} =\bigcap_{l\in\Delta}l\mathbb{Z}_{(l)}.$$
In particular, if $\Delta\neq\emptyset$
then $J_{a,b}+J_{c,d}$ is the Jacobson radical $J(R)$ of the semilocal ring $R$. 
\end{lemma}
{\em Proof:}
Let $\Delta_{a,b}^\prime:=\left\{\begin{array}{ll}
\Delta_{a,b}\setminus\{ 2,\infty\} & \mbox{if }2\in\Delta_{a,b}\mbox{ and }
v_2(a), v_2(b)\mbox{ are even}\\
\Delta_{a,b}\setminus\{\infty\} & \mbox{else}
\end{array}\right.$, and similarly $\Delta_{c,d}^\prime$.
Then, by Lemma \ref{Jab}(d) (for the first equality)
and by weak approximation (for the second),
$$J_{a,b}+J_{c,d} = \bigcap_{l\in\Delta_{a,b}^\prime}l\mathbb{Z}_{(l)}
+  \bigcap_{l\in\Delta_{c,d}^\prime}l\mathbb{Z}_{(l)}
= \bigcap_{l\in\Delta_{a,b}^\prime\cap\Delta_{c,d}^\prime}l\mathbb{Z}_{(l)}.$$
By our assumption on $a,b,c,d$, however,
$\Delta_{a,b}\cap\Delta_{c,d}=\Delta_{a,b}^\prime\cap\Delta_{c,d}^\prime$,
which proves the first claim.

The `in particular' follows immediately.
\qed

Now let us first turn to the rings $R_p^{[k]}$ for $k=3,5$ and $7$
defined in Definition \ref{R} and recall that 
$$R_p^{[k]} =\left\{\begin{array}{ll}
T_{-1,-p}+T_{2,-p} & \mbox{if }k=3\\
T_{-2,-p}+T_{2,-p} & \mbox{if }k=5\\
T_{-1,-p}+T_{-2,p} & \mbox{if }k=7
\end{array}\right.$$
\begin{corollary}
\label{Phi_k}
Define for $k=1,3,5$ and $7$,
$$\begin{array}{ccl}
\Phi_k & := & \{p\in {\mathbb Q}^{>0}\mid p\equiv k\;(\!\!\!\mod 8\mathbb{Z}_{(2)})\mbox{ and } {\mathbb P}(p)\subseteq {\mathbb P}^{[1]}\cup {\mathbb P}^{[k]}\}\\
\Psi & := & \{(p,q)\in \Phi_1\times \Phi_3\mid
p\in 2\cdot ({\mathbb Q}^\times)^2\cdot (1+ J(R_q^{[3]}))\}.
\end{array}$$

\begin{enumerate}
\item[(a)]
Then $\Phi_k$ is diophantine in ${\mathbb Q}$.
\item[(b)]
If $k=3,5$ or $7$ and if $p\in\Phi_k$ then
$\mathbb{P}^{[k]}(p)\neq\emptyset$ and
$$\{ 0\}\neq J(R_p^{[k]})=\left\{\begin{array}{ll}
J_{-1,-p}+J_{2,-p} & \mbox{if }k=3\\
J_{-2,-p}+J_{2,-p} & \mbox{if }k=5\\
J_{-1,-p}+J_{-2,p} & \mbox{if }k=7
\end{array}\right.$$
In particular, in each of the cases,
the Jacobson radical is diophantine in $\mathbb{Q}$,
by a formula that is uniform in $p$.
\item[(c)]
$\Psi$ is diophantine in ${\mathbb Q}$.
\end{enumerate}
\end{corollary}
{\em Proof:}
(a) It is clear that `$p>0$' is diophantine.
It is also clear from Proposition \ref{R_p}(a) that, for $k=1,3,5$ and $7$, the property `$p\equiv k\;(\!\!\!\mod 8\mathbb{Z}_{(2)})$' is diophantine.

Moreover, if $v_2(p)$ is even and $k^\prime =3,5$ or $7$, then, by Proposition \ref{R_p}(b),
$${\mathbb P}^{[k^\prime]}(p) =\emptyset \Longleftrightarrow p\in ({\mathbb Q}^\times)^2\cdot (R_p^{[k^\prime]})^\times$$
(Note that we are {\em not} assuming that $p\equiv k^\prime \;(\!\!\!\mod 8\mathbb{Z}_{(2)}$.)
So the property on the left is diophantine. But then so are
$$\begin{array}{l}
\Phi_1 =\{p\equiv 1\;(\!\!\!\mod 8\mathbb{Z}_{(2)})\mid {\mathbb P}_3(p)=\emptyset,\,{\mathbb P}_5(p)=\emptyset\mbox{ and }{\mathbb P}_7(p)=\emptyset\}\\
\Phi_3 =\{p\equiv 3\;(\!\!\!\mod 8\mathbb{Z}_{(2)})\mid {\mathbb P}^5(p)=\emptyset\mbox{ and }{\mathbb P}^7(p)=\emptyset\}\\
\Phi_5 =\{p\equiv 5\;(\!\!\!\mod 8\mathbb{Z}_{(2)})\mid {\mathbb P}^3(p)=\emptyset\mbox{ and }{\mathbb P}^7(p)=\emptyset\}\\
\Phi_7 =\{p\equiv 7\;(\!\!\!\mod 8\mathbb{Z}_{(2)})\mid {\mathbb P}^3(p)=\emptyset\mbox{ and }{\mathbb P}^5(p)=\emptyset\}.
\end{array}$$
(b) Assume $k=3,5$ or $7$ and that $p\in\Phi_k$.
Then $p\equiv k\;(\!\!\!\mod 8\mathbb{Z}_{(2)})$
and so, by Proposition \ref{R_p}(b), $R_p^{[k]}=\bigcap_{l\in\mathbb{P}^{[k]}(p)}\mathbb{Z}_{(l)}$.
As $p>0$ and $p\equiv k\;(\!\!\!\mod 8\mathbb{Z}_{(2)})$, $\mathbb{P}^{[k]}(p)\neq\emptyset$
and hence $J(R_p^{[k]})=\bigcap_{l\in\mathbb{P}^{[k]}(p)}l\mathbb{Z}_{(l)}\neq\{ 0\}$. The explicit formulas now follow from Lemma \ref{J+J}, as the assumptions of the Lemma are satisfied in each case.
(c) follows directly from (a) and (b).\qed
\medskip

The most difficult case is when $p\in\Phi_1$.
Recall from Definition \ref{R} and from Proposition \ref{R_p}(c) that, for $p,q\in\mathbb{Q}^\times$, we have defined
$R_{p,q}^{[1]}:= T_{-2p,q}+T_{2p,q}$
and $\mathbb{P}(p,q):=\Delta_{-2p,q}\cap\Delta_{2p,q}$.
\begin{proposition}
\label{RnotQ}
\begin{enumerate}
\item[(a)]
If $(p,q)\in\Psi$, then $\mathbb{P}(p,q)\neq\emptyset$.
\item[(b)]
If $(p,q)\in\Psi$, then $J(R_{p,q}^{[1]})=J_{-2p,q}+J_{2p,q}$.
\item[(c)]
The set $\{ (p,q,x)\in\mathbb{Q}^3\mid (p,q)\in\Psi\mbox{ and }x\in J(R_{p,q}^{[1]})\}$ is diophantine.
\end{enumerate}
\end{proposition}
{\em Proof:}
(a) Assume $(p,q)\in\Psi$.
Multiplying $p$ or $q$ by nonzero rational squares does not change
$R_{p,q}^{[1]}$ or $J_{-2p,q}$ or $J_{2p,q}$,
so we can assume that $p$ and $q$ are squarefree positive integers.
Since $p\equiv 1\;(\!\!\!\mod 8\mathbb{Z}_{(2)})$
and $q\equiv 3\;(\!\!\!\mod 8\mathbb{Z}_{(2)}$,
we have, by Observation \ref{pinDelta}, $(2p,q)_2=-1$.
By Hilbert reciprocity, there must also be an odd prime $l$
such that $(2p,q)_l=-1$.
By definition of $\Psi$ and, again, by Observation \ref{pinDelta},
this implies that $l\in\{ 1,3\} + 8\mathbb{Z}_{(2)}$
and $l\not\in\mathbb{P}^{[3]}(q)$.
These two conditions imply $(-1,q)_l=1$.
Multiplying yields $(-2p,q)_l=-1$.
Thus $l\in\mathbb{P}(p,q)$.\\
(b) is immediate from (a) and Lemma \ref{J+J}.\\
(c) follows from Corollary \ref{Phi_k}(c), from (b) and the note preceding Lemma \ref{Jab}.\qed 
\subsection*{{\em Step 4:} From existential to universal}
Let $R$ be a semilocal subring of ${\mathbb Q}$, i.e., $R=\bigcap_{\l\in\Delta}{\mathbb Z}_{(l)}$ for some finite $\Delta\subseteq {\mathbb P}$. Define
$$\widetilde{R}:=\{x\in {\mathbb Q}\mid\neg\exists y\in J(R)\mbox{ with }x\cdot y=1\}.$$
\begin{lemma}
\label{widetildeR}
\begin{enumerate}
\item[(a)]
If $J(R)$ is diophantine in ${\mathbb Q}$ then $\widetilde{R}$ is defined by a {\em universal} formula in ${\mathbb Q}$.
\item[(b)]
$\widetilde{R}=\bigcup_{l\in\Delta}{\mathbb Z}_{(l)}$,
provided $\Delta\neq\emptyset$, i.e., provided $R\neq {\mathbb Q}$.
\item[(c)]
In particular, if $R={\mathbb Z}_{(p)}$ for some $p\in {\mathbb P}$ then $\widetilde{R}=R$.
\end{enumerate}
\end{lemma}
{\em Proof:}
(a) is obvious from the definition of $\widetilde{R}$, and (c) is a special case of (b). So we only need to prove (b).

For the inclusion `$\subseteq$', pick $x\in \widetilde{R}$ and assume that $x\not\in \bigcup_{l\in\Delta}{\mathbb Z}_{(l)}$.
Then for all $l\in\Delta$, $v_l(x)<0$, and hence $y:=x^{-1}\in\bigcap_{l\in\Delta}l\,{\mathbb Z}_{(l)} = J(R)$, contradicting our assumption that $x\in\widetilde{R}$.

For the converse inclusion `$\supseteq$', assume $x\in {\mathbb Z}_{(l)}$ for some $l\in\Delta$. Then, for any $y\in J(R)$, $x\cdot y\in l\,{\mathbb Z}_{(l)}$, so, in particular $x\cdot y\neq 1$.\qed
\medskip

Now we can give our universal definition of ${\mathbb Z}$ in ${\mathbb Q}$:
\begin{proposition}
\begin{enumerate}
\item[(a)] 
$${\mathbb Z}=\widetilde{{\mathbb Z}_{(2)}}\cap\left( \bigcap_{k=3,5,7}\bigcap_{p\in\Phi_k}\widetilde{R_p^{[k]}}\right)
\cap\bigcap_{(p,q)\in\Psi}\widetilde{R_{p,q}^{[1]}},$$
where $\Phi_k$ and $\Psi$ are the diophantine sets defined in Corollary \ref{Phi_k}.
\item[(b)]for any $t\in {\mathbb Q}$,
$$\begin{array}{lll}t\in{\mathbb Z} & \Longleftrightarrow & t\in\widetilde{{\mathbb Z}_{(2)}}\wedge\\
 & &\forall p\bigwedge_{k=3,5,7} (t\in\widetilde{R_p^{[k]}}\vee p\not\in\Phi_k)\wedge\\
 & &\forall p,q (t\in \widetilde{R_{p,q}^{[1]}}\vee (p,q)\not\in\Psi)
\end{array}$$
\item[(c)] ({\bf Theorem 1})
There is a natural number $n$ and a polynomial $g\in {\mathbb Z}[t;x_1,\ldots ,x_n]$ such that, for any $t\in {\mathbb Q}$,
$$t\in {\mathbb Z}\mbox{ iff }\forall x_1\ldots\forall x_n\in {\mathbb Q}\; g(t;x_1,\ldots ,x_n)\neq 0.$$
\end{enumerate}
\end{proposition}
{\em Proof:}
(a) The equation is valid by Proposition \ref{R_p} and Lemma \ref{widetildeR}(b), (c).

(b) This is a reformulation of (a)
revealing that the formula thus obtained for ${\mathbb Z}$ in ${\mathbb Q}$
{\em is} universal:
the $\widetilde{R}$'s are universal by Corollary \ref{Phi_k}, Proposition \ref{RnotQ} and Lemma \ref{widetildeR}(a); $\Phi_k$ and $\Psi$ are existential by Corollary \ref{Phi_k}(a) and (c), so their negation is universal as well.

(c) This is immediate from (b).\qed
\section{More diophantine predicates in ${\mathbb Q}$}
From the results and techniques of section 2, one obtains new diophantine predicates in ${\mathbb Q}$.
They are of interest in their own right, but maybe they can also be used to show that Hilbert's 10th problem over $\mathbb{Q}$ cannot be solved, not by defining or interpreting ${\mathbb Z}$ in ${\mathbb Q}$, but, e.g., by assigning graphs to the various finite sets of primes encoded in these predicates, and using graph theoretic undecidability results.
We will also use some of these new predicates for our $\forall\exists$-definition of ${\mathbb Z}$ in ${\mathbb Q}$ which uses just one universal quantifier (Corollary \ref{just1}).

Before listing the new diophantine predicates we shall first prove the following
\begin{lemma}
\label{R_1}
Assume $p\in\Phi_1$ and define\footnote{We hope the notation $R_p^{[1]}$
is not too confusing as the definition is different from that of $R_p^{[k]}$ for $k=3,5$ or $7$. The crucial property, however, the `in particualr', gives the same result as in Corollary \ref{R_p}(b).}
$$R_p^{[1]}:=\{ x\in {\mathbb Q}\mid\exists q\mbox{ with }(p,q)\in\Psi, q\in (R_{p,q}^{[1]})^\times\mbox{ and }x\in R_{p,q}^{[1]}\}.$$
Then $R_p^{[1]}$ is diophantine in ${\mathbb Q}$
and $R_p^{[1]}=\bigcup_{l\in {\mathbb P}(p)} {\mathbb Z}_{(l)}$
(which is $\emptyset$ if $\mathbb{P}(p)=\emptyset$, i.e., if $p\in\mathbb{Q}^2$).

In particular, if $p$ is a prime $\equiv 1\!\!\mod 8$ then $R_p^{[1]}={\mathbb Z}_{(p)}$.
\end{lemma}
{\em Proof:}
That $R_p^{[1]}$ is diophantine in ${\mathbb Q}$ is immediate from
Corollary \ref{Phi_k}.

Assuming $(p,q)\in\Psi$, the condition `$q\in (R_{p,q}^{[1]})^\times$' implies that,
in the terminology of Proposition \ref{R_p}(c),
${\mathbb P}(p,q)\subseteq {\mathbb P}(p)$:
Suppose $q\in (R_{p,q}^{[1]})^\times$ and $l\in\mathbb{P}(p,q)$.
Then $v_l(q)=0$, so $l\not\in\mathbb{P}(q)$.
Hence, by the definition of $\mathbb{P}(p,q)$,
$l\in\mathbb{P}(p)$ (note that, by Observation \ref{pinDelta}, $2\not\in\mathbb{P}(p,q)$, as $(p,q)\in\Phi_1\times\Phi_3$).

This yields the last inclusion in
$$R_p^{[1]}=\bigcup_{\begin{array}{ll}
q\in (R_{p,q}^{[1]})^\times\mbox{\small with}\\
(p,q)\in\Psi
\end{array}}
R_{p,q}^{[1]}\;\;= \bigcup_{\begin{array}{ll}
q\in (R_{p,q}^{[1]})^\times\mbox{\small with}\\
(p,q)\in\Psi
\end{array}}
\bigcap_{l\in {\mathbb P}(p,q)}{\mathbb Z}_{(l)}
\;\;\subseteq\bigcup_{l\in {\mathbb P}(p)}{\mathbb Z}_{(l)}.$$
The first equality is by definition, the second by Proposition \ref{R_p}(c)
using that, from Proposition \ref{RnotQ}(a),  ${\mathbb P}(p,q)\neq\emptyset$.

Conversely, suppose $l\in {\mathbb P}(p)$ and $x\in {\mathbb Z}_{(l)}$.
Choose a prime $q\equiv 3\!\!\mod 8$ with $\leg{l}{q}=-1$
and with $\leg{l^\prime}{q}=1$
for each $l^\prime\in {\mathbb P}(p)\setminus\{ l\}$.

Then $\leg{pq^{-v_q(p)}}{q} = -1$, so $(p,q)\in\Psi$:
note that $v_q(p)$ is even since $p\in\Phi_1$,
so $\phi_q(pq^{-v_q(p)})$ is a non-square in $\mathbb{F}_q$,
i.e., $\in 2\cdot (\mathbb{F}_q^\times)^2$;
hence $p\in 2\cdot (\mathbb{Q}^\times)^2(1+q\mathbb{Z}_{(q)})$.

Therefore, by Proposition \ref{R_p}(c)
and the Quadratic Reciprocity Law,
${\mathbb P}(p,q) = \{ l\}$:
Clearly $l\in\mathbb{P}(p,q)$ as $l\in\mathbb{P}(p)$
and $\leg{q}{l}=\leg{l}{q}=-1$ ($p\in\Phi_1$, so $l\equiv 1\!\!\mod 8$);
and for any $l^\prime\in\mathbb{P}(p)\setminus \{ l\}$,
$\leg{l^\prime}{q}=\leg{q}{l^\prime}=1$, so
$\l^\prime\not\in\mathbb{P}(p,q)$;
finally $\mathbb{P}(q)=\{ q\}$, but $\leg{2pq^{-v_q(2p)}}{q}=1$,
hence $q\not\in\mathbb{P}(p,q)$.

Thus $x\in\mathbb{Z}_{(l)} = R_{p,q}^{[1]}$ and 
$q\in (R_{p,q}^{[1]})^\times$.\qed

\begin{proposition}
\label{new}
For $x,y\in {\mathbb Q}^\times$, the following properties are diophantine:
\begin{enumerate}
\item[(a)]
for fixed $k\in\{ 3,5,7\}$, the property that
$x,y\in\Phi_k$ and ${\mathbb P}^{[k]}(x)\cap {\mathbb P}^{[k]}(y)=\emptyset$
\item[(b)] $x\not\in {\mathbb Q}^2$
\item[(c)] for fixed $k\in\{ 1,3,5,7\}$, the property that
$x\equiv k\;(\!\!\!\mod 8\mathbb{Z}_{(2)})$ and $x\not\in\Phi_k$
\item[(d)] for fixed $k\in\{ 3,5,7\}$, the property that
${\mathbb P}^{[k]}(x) =\emptyset$
\item[(e)] $x\not\in N(y)$, where $N(y)$ is the image of the norm ${\mathbb Q}(\sqrt{y})\to{\mathbb Q}$
\end{enumerate}
\end{proposition}
{\em Proof:}
{\bf (a)}
By Corollary \ref{Phi_k}(a), $\Phi_k$ is diophantine.
By Corollary \ref{Phi_k}(b), for any $x\in\Phi_k$, ${\mathbb P}^{[k]}(x)\neq\emptyset$ and hence 
$J(R_x^{[k]})$ is diophantine.
Now let $x,y\in\Phi_k$ and
recall that, by Proposition \ref{R_p},
$R_x^{[k]}=\bigcap_{l\in\mathbb{P}^{[k]}(x)}\mathbb{Z}_{(l)}$,
and likewise for $R_y^{[k]}$.
So we have the equivalence
$${\mathbb P}^{[k]}(x)\cap {\mathbb P}^{[k]}(y)=\emptyset\Longleftrightarrow
1\in J(R_x^{[k]})+J(R_y^{[k]})$$
which then proves {\bf (a)}.
\medskip

{\bf (b)}
The property that `$v_2(x)$ is odd' is diophantine:
$v_2(x)$ is odd if and only if $x=2yz^2$ 
for some $y\in\mathbb{Z}_{(2)}^\times$ and some $z\in\mathbb{Q}^\times$.
As the property `$x<0$' is diophantine as well,
by Corollary \ref{Phi_k}(a) and (b), it suffices to show
$$x\not\in {\mathbb Q}^2\Longleftrightarrow\left\{\begin{array}{l}
x<0\mbox{ or }v_2(x)\mbox{ is odd or}\\
\exists p\in\Phi_3\mbox{ with }x\in 2\cdot ({\mathbb Q}^\times)^2\cdot (1+J(R_p^{[3]}))
\end{array}\right.$$
`$\Rightarrow$':
Assume that $x\not\in {\mathbb Q}^2$, that $x>0$ and that $v_2(x)$ is even.
Multiplying $x$ by a nonzero rational square does not change
the truth of either side of the implication, so we may assume that
$x=p_1\cdots p_r$ for distinct odd primes $p_1,\ldots ,p_r$
where $r\geq 1$.
$$\mbox{Choose }a_1\in {\mathbb Z}\mbox{ with }\leg{a_1}{p_1} =
\left\{\begin{array}{ll}
-1 & \mbox{if }p_1\equiv 1\!\!\mod 4\\
1 & \mbox{if }p_1\equiv 3\!\!\mod 4
\end{array}\right.$$
and, for $i>1$,
$$\mbox{choose }a_i\in {\mathbb Z}\mbox{ with }\leg{a_i}{p_i} =
\left\{\begin{array}{ll}
1 & \mbox{if }p_i\equiv 1\!\!\mod 4\\
-1 & \mbox{if }p_i\equiv 3\!\!\mod 4
\end{array}\right.$$
Finally, choose a prime $p\equiv 3\!\!\mod 8$ with $p\equiv a_i\!\!\mod p_i$ ($i=1,\ldots ,r$).
Then, by the Quadratic Reciprocity Law, $\leg{x}{p}=-1$.

Clearly, $p\in\Phi_3$. By Lemma \ref{R_p}(b), $R_p^{[3]}={\mathbb Z}_{(p)}$.
Hence $x\in 2\cdot ({\mathbb Q}^\times)^2\cdot (1+J(R_p^{[3]}))$,
as $\leg{2}{p}=-1$.

`$\Leftarrow$':
If $x<0$ or $v_2(x)$ is odd then clearly $x\not\in {\mathbb Q}^2$.

If $x\in 2\cdot ({\mathbb Q}^\times)^2\cdot (1+J(R_p^{[3]}))$
for some $p\in\Phi_3$ then ${\mathbb P}^{[3]}(p)\neq\emptyset$,
and for any $l\in {\mathbb P}^{[3]}(p)$ one has
that $v_l(x)$ is even and
$\leg{xl^{-v_l(x)}}{l} = \leg{2}{l} =-1$.
Hence $x\not\in {\mathbb Q}^2$.
\medskip

{\bf (c)}
By Proposition \ref{R_p}(a), $x\equiv k\;(\!\!\!\mod 8\mathbb{Z}_{(2)})$ is diophantine.
First assume $x\equiv 1\;(\!\!\!\mod 8\mathbb{Z}_{(2)})$.
Then $x\not\in\Phi_1$ if and only if
$x\leq 0$ or $x>0$ and, for some $k\in\{ 3,5,7\}$, $\mathbb{P}^{[k]}(x)\neq\emptyset$.

This last condition can be expressed diophantinely
by distinguishing the cases whether the number of $k\in\{ 3,5,7\}$
with $\mathbb{P}^{[k]}(x)\neq\emptyset$ is 1, 2 or 3.

If it is 1, say $\mathbb{P}^{[k]}(x)\neq\emptyset$,
then $\sharp\mathbb{P}^{[k]}(x)$ must be even
(in order to get $x\equiv 1 \;(\!\!\!\mod 8\mathbb{Z}_{(2)})$),
so we can choose $p\in\mathbb{P}^{[k]}(x)$
and let 
$$ y_k:=p^{v_p(x)}\mbox{ and }
y_k^\prime:=\prod_{l\in\mathbb{P}^{[1]}(x)}l^{v_l(x)} \cdot
\prod_{l\in\mathbb{P}^{[k]}(x)\setminus \{p\}}l^{v_l(x)}.$$
Then $y_k, y_k^\prime\in\Phi_k$, $\mathbb{P}^{[k]}(y_k)\cap\mathbb{P}^{[k]}(y_k^\prime)=\emptyset$
and $x=y_k\cdot y_k^\prime$.
By (a), the condition that there exist such $y_k, y_k^\prime$ is diophantine,
and, when satisfied, it implies $x\not\in\Phi_1$.

If $\{ k\in\{ 3,5,7\}\mid\mathbb{P}^{[k]}(x)\neq\emptyset\} =\{ k_1,k_2\}$ 
for distinct $k_1,k_2$ then
both $\sharp\mathbb{P}^{[k_1]}(x)$ and $\sharp\mathbb{P}^{[k_2]}(x)$
must be even, again, and so one constructs similarly $y_1,y_1^\prime\in\Phi_{k_1}$ and $y_2,y_2^\prime\in \Phi_{k_2}$ with $\mathbb{P}^{[k_i]}(y_i)\cap\mathbb{P}^{[k_i]}(y_i^\prime)=\emptyset$
for $i=1,2$ such that $x=y_1\cdot y_1^\prime\cdot y_2\cdot y_2^\prime$.

If $\mathbb{P}^{[k]}(x)\neq\emptyset$ for all three $k\in \{ 3,5,7\}$
then either all three sets have an even number of elements or all three have an odd number of elements, and in either case it is clear how to proceed along the same lines.

Now assume $x\equiv 3\;(\!\!\!\mod 8\mathbb{Z}_{(2)})$.
Then $x\not\in\Phi_3$ if and only if
$x\leq 0$ or $x>0$ and $\mathbb{P}^{[5]}(x)\neq\emptyset$
or $\mathbb{P}^{[7]}(x)\neq\emptyset$.
Here the last condition is diophantine again, distinguishing the cases whether the number of $k\in\{ 5,7\}$ with $\mathbb{P}^{[k]}(x)\neq\emptyset$
is 1 or 2 etc.

It is clear how similar existential formulas can be written down for `$x\not\in\Phi_5$' and `$x\not\in\Phi_7$'.
\medskip

{\bf (d)}
${\mathbb P}^{[3]}(x)=\emptyset$ if and only if,
modulo a nonzero rational square factor,
$x$ or $-x$ or $2x$ or $-2x$ is a product of primes in $\bigcup_{k=1,5,7}\mathbb{P}^{[k]}$.
Note that for a fixed $k\in \{1,5,7\}$,
each product of primes in $\mathbb{P}^{[k]}$
can be expressed as a product of one or two factors of elements in $\Phi_k$.
Hence ${\mathbb P}^{[3]}(x)=\emptyset$ if and only if
$$\begin{array}{c}
\exists y_1,\ldots ,y_8,z\\
\left( y_1,y_2\in\Phi_1\wedge y_3,y_4\in\Phi_5\wedge y_5,y_6\in\Phi_7\wedge y_7=-1\wedge y_8=2\right.\\
\left.\wedge\bigvee_{I\subseteq \{1,\ldots ,8\}} x=z^2\prod_{i\in I}y_i\right)
\end{array}$$
And, again, similar formulas hold for $k=5$ and $k=7$.
\medskip

{\bf (e)}
$x\not\in N(y)$ iff
$$\begin{array}{l}
(x<0\wedge y<0)\\
\vee\bigvee_{k=3,5,7}\exists p\in\Phi_k\mbox{ with}\\
\left(\left(x\in p\cdot ({\mathbb Q}^\times)^2\cdot (R_p^{[k]})^\times\right)\wedge
\left( y\mbox{ or }-xy\in a_k\cdot ({\mathbb Q}^\times )^2\cdot (1+J(R_p^{[k]}))\right)\right.\\
\left. \vee\left( y\in p\cdot ({\mathbb Q}^\times)^2\cdot (R_p^{[k]})^\times\right)\wedge
\left( x\mbox{ or }-xy\in a_k\cdot ({\mathbb Q}^\times )^2\cdot (1+J(R_p^{[k]}))\right)\right)\\
\vee\exists (p,q)\in\Psi\mbox{ with }q\in (R_{p,q}^{[1]})^\times\mbox{ and}\\
\left(\left(x\in p\cdot ({\mathbb Q}^\times)^2\cdot (R_{p,q}^{[1]})^\times\right)\wedge
\left( y\mbox{ or }-xy\in q\cdot ({\mathbb Q}^\times )^2\cdot (1+J(R_{p,q}^{[1]}))\right)\right.\\
\left. \vee\left( y\in p\cdot ({\mathbb Q}^\times)^2\cdot (R_{p,q}^{[1]})^\times\right)\wedge
\left( x\mbox{ or }-xy\in q\cdot ({\mathbb Q}^\times )^2\cdot (1+J(R_{p,q}^{[1]}))\right)\right)
\end{array}$$
where $a_3=a_5=2$ and $a_7=-1$.

This uses Observation \ref{pinDelta}(b) and (c), Corollary \ref{Phi_k}(b) and (c), the previous parts and the local-global principle for norms. 

The first line says that $x\not\in N(y)$ over ${\mathbb R}$.

Lines 2-4 say that $x\not\in N(y)$ over ${\mathbb Q}_l$
for some non-empty set of primes $l\equiv 3,5$ or $7\!\!\mod 8$:
Fix $k\in \{ 3,5,7\}$.
By Corollary \ref{Phi_k}(b), $p\in\Phi_k$ implies that $\mathbb{P}^{[k]}(p)\neq\emptyset$.
We claim that 
$$ (x,y)_l=-1\mbox{ for some }l\in\mathbb{P}^{[k]}
\Longleftrightarrow\exists p\in\Phi_k\mbox{ with }\left(\ldots\right),$$
where $(\ldots )$ is the bracket is line 3 and 4.\\
`$\Rightarrow$':
Assume $l\in\mathbb{P}^{[k]}$ with $(x,y)_l=-1$.
Let $p=l$. Then $R_p^{[k]}=\mathbb{Z}_l$ and `$(\ldots)$'
says that $v_l(x)$ is odd and $yl^{-v_l(y)}$ or $-xyl^{-v_l(xy)}$
is a quadratic non-residue $\!\!\mod l$ or the same with $x$ and $y$ swapped. By Observation \ref{pinDelta}, this is equivalent to $(x,y)_l=-1$,
so it holds by our assumption.\\
`$\Leftarrow$':
Suppose $p\in\Phi_k$ satisfies `$(\ldots )$'.
Then $\mathbb{P}^{[k]}(p)\neq\emptyset$
and, for any $l\in\mathbb{P}^{[k]}(p)$,
$v_l(x)$ is odd and, by the choice of $a_k$, either $yl^{-v_l(y)}$ or $-xyl^{-v_l(xy)}$
is a quadratic non-residues $\!\!\mod l$ or the same with $x$ and $y$ swapped,
so $(x,y)_l=-1$.

Lines 5-7 say that $x\not\in N(y)$ over ${\mathbb Q}_l$
for some non-empty set of primes $l\equiv 1\!\!\mod 8$.
As in the proof of Lemma \ref{R_1},
the condition `$q\in (R_{p,q}^{[1]})^\times$'
makes sure that, in the terminology of Proposition \ref{R_p}(c),
${\mathbb P}(p,q)\cap {\mathbb P}(q) = \emptyset$,
so ${\mathbb P}(p,q)\subseteq {\mathbb P}(p)$.
And, by Proposition \ref{RnotQ}(a), ${\mathbb P}(p,q)\neq\emptyset$.
Line 6 and 7 then say that $x\not\in N(y)$ over ${\mathbb Q}_l$
for any $l\in {\mathbb P}(p,q)$.
Note that the role of $a_k$ in lines 3 and 4 of being a quadratic non-residue $\!\!\mod l$ for all $l\in\mathbb{P}^{[k]}$
is here taken by $q$
which is a quadratic non-residue for all $l\in\mathbb{P}^{[1]}(p)$
with $(p,q)\in\Psi$.

We could disregard the prime $p=2$, as `$x\not\in N(y)$' either happens nowhere locally, or at least at two primes in ${\mathbb P}\cup\{\infty\}$.\qed
\medskip

The result in {\bf (b)} was also obtained in [P2] -- using a deep result on Ch\^atelet surfaces from [CSS] -- our proof is elementary.

Let us also mention that {\bf (b)} follows from {\bf (e)}:
$x\not\in {\mathbb Q}^2\Leftrightarrow\exists y\; x\not\in N(y)$
(and we did not use {\bf (b)} in order to prove {\bf (e)}).
\medskip

We close this section by showing that there is an $\forall\exists$-definition of ${\mathbb Z}$ in ${\mathbb Q}$ with just one universal quantifier:
\begin{corollary}
\label{just1}
For all $t\in {\mathbb Q}$, $t\in{\mathbb Z}$ if and only if 
$$\forall p \left( t\in {\mathbb Z}_{(2)}\,\wedge\,
\left\{
\begin{array}{l}
\;\;\;\left( p\in {\mathbb Q}^2\cdot (2 + 4 {\mathbb Z}_{(2)})\right)\\
\vee\bigvee_{k=1,3,5,7} \left\{
\begin{array}{l} \;\;\;\left( p\neq 0\wedge p\in{\mathbb Q}^2\cdot (k+8{\mathbb Z}_{(2)})\right)\\
\wedge \left( \left( p\not\in\Phi_k\right)\vee p\in\mathbb{Q}^2\vee \left( p\in\Phi_k\setminus {\mathbb Q}^2 \wedge t\in R_p^{[k]}\right)\right)\end{array}\right.
\end{array}\right.\right)$$
\end{corollary}
{\em Proof:}
The equivalence holds by Proposition \ref{R_p}(a) and (b) and by Lemma \ref{R_1}.

That the resulting formula is of the shape $\forall\exists$ with just one universal quantifier `$\forall p$' follows from Proposition \ref{R_p}, Corollary \ref{Phi_k}, Lemma \ref{R_1} and Proposition \ref{new}.
Note that, under the assumption `$p\in{\mathbb Q}^2\cdot (k+8{\mathbb Z}_{(2)})$', the property `$p\not\in\Phi_k$'
is equivalent to `$p\not\in\mathbb{Z}_{(2)}^\times$ or
($p\in k+\mathbb{Z}_{(2)}\mbox{ and }p\not\in\Phi_k$)'
which is diophantine by Proposition \ref{new}(c).
And `$p\not\in\mathbb{Q}^2$' is diophantine by \ref{new}(b).\qed
\section{Why ${\mathbb Z}$ should not be diophantine in ${\mathbb Q}$}
In this section we show that ${\mathbb Z}$ is not diophantine in ${\mathbb Q}$, provided one believes in a certain version of what one may (arguably) call `the Bombieri-Lang Conjecture' on varieties with many rational points.

The version of this conjecture in the special case of varieties over ${\mathbb Q}$ on which our result is based is the following (mainly after section F.5.2 of [HS]):
\medskip\\
{\bf Bombieri-Lang Conjecture}
{\em Let $V$ be an absolutely irreducible affine or projective positive-dimensional variety over ${\mathbb Q}$
such that $V({\mathbb Q})$ is Zariski dense in $V$.
Then so is
$$\bigcup_{\phi:A\; \dashrightarrow V} \phi (A({\mathbb Q})),$$
where the $\phi:A\; \dashrightarrow V$ run through all non-constant ${\mathbb Q}$-rational maps from positive-dimensional abelian varieties $A$ defined over ${\mathbb Q}$ to $V$.}

\begin{lemma}
\label{Bombieri-Lang}
Assume the Bombieri-Lang Conjecture as above.
Let $f\in\mathbb{Q}[x_1,\ldots ,x_{n+1}]\setminus\mathbb{Q}[x_1,\ldots ,x_n]$
be absolutely irreducible and let
Let $V=V(f)\subseteq {\mathbb A}^{n+1}$ be the affine hypersurface defined over ${\mathbb Q}$ by $f$.
Assume that $V({\mathbb Q})$ is Zariski dense in $V$.
Let $\pi:{\mathbb A}^{n+1}\to {\mathbb A}^1$ be the projection onto the first coordinate.
Then $V({\mathbb Q})\cap\pi^{-1}({\mathbb Q}\setminus {\mathbb Z})$ is also Zariski dense in $V$.
\end{lemma}
(For $n=1$ the Lemma holds unconditionally, by Siegel's Theorem.)
\medskip\\
{\em Proof:}
Choose any $g\in{\mathbb Q}[x_1,\ldots ,x_n]\setminus \{ 0\}$.
By the Bombieri-Lang Conjecture there are an abelian variety $A$
and a rational map $\phi:\; A\dashrightarrow  V$, both defined over ${\mathbb Q}$,
such that $\phi(A({\mathbb Q}))\setminus V(g)({\mathbb Q})$ is infinite
(considering $V(g)$ as subset of ${\mathbb A}^{n+1}$).
By possibly composing $\phi$ with another rational map (from the left),
we may assume that
$\pi (\phi (A({\mathbb Q}))\setminus V(g)({\mathbb Q}))$ is infinite,
and that the pole divisor $D$ of $\pi\circ\phi$ is ample.
By Corollary 6.2 in [F], there are only finitely many $P\in A({\mathbb Q})\setminus D({\mathbb Q})$ with $\pi (\phi (P))\in {\mathbb Z}$
(cf. the remarks following Theorem 1 in [Si]).
This implies that $(V({\mathbb Q})\setminus V(g)({\mathbb Q}))\cap\pi^{-1}({\mathbb Q}\setminus {\mathbb Z})\neq\emptyset$.
Since $g$ was arbitrary this shows that
$V({\mathbb Q})\cap\pi^{-1}({\mathbb Q}\setminus {\mathbb Z})$ is Zariski dense in $V$.\qed

\begin{corollary}
\label{notex}
Assume the Bombieri-Lang Conjecture as stated above.
Then there is no infinite subset of ${\mathbb Z}$ existentially definable in ${\mathbb Q}$.
In particular, ${\mathbb Z}$ is not diophantine in ${\mathbb Q}$.
\end{corollary}
{\em Proof:}
If $\mathbb{Z}$ contains an infinite subset
that is diophantine over $\mathbb{Q}$
then there is a hypersurface $W$ in $\mathbb{A}^{n+1}$
such that $\pi (W(\mathbb{Q}))$ is infinite
(where $\pi$ is as in Lemma \ref{Bombieri-Lang}).
Replace $W$ by the Zariski closure $\overline{W}$ of $W(\mathbb{Q})$;
this ensures that the irreducible components $V$ of $\overline{W}$ are geometrically irreducible (given by absolutely irreducible polynomials).
For at least one such $V$ the set $\pi (V(\mathbb{Q}))$ is infinite
(and still contained in $\mathbb{Z}$:
note that $W(\mathbb{Q})=\overline{W}(\mathbb{Q})$).
This contradicts Lemma \ref{Bombieri-Lang}.
\qed
\medskip

Let us conclude with a collection of closure properties for pairs of models of Th$({\mathbb Q})$ (in the ring language), one a substructure of the other,
which might have a bearing on the final (unconditional) answer
to the question whether or not ${\mathbb Z}$ is diophantine in ${\mathbb Q}$.
\begin{proposition}
Let ${\mathbb Q}^\star, {\mathbb Q}^{\star\star}$ be models of Th$({\mathbb Q})$ (i.e. elementary extensions of ${\mathbb Q}$) with ${\mathbb Q}^\star\subseteq {\mathbb Q}^{\star\star}$, and let ${\mathbb Z}^\star$ and ${\mathbb Z}^{\star\star}$ be their rings of integers. Then
\begin{enumerate}
\item[(a)]
${\mathbb Z}^{\star\star}\cap {\mathbb Q}^\star\subseteq {\mathbb Z}^\star$.
\item[(b)]
${\mathbb Z}^{\star\star}\cap {\mathbb Q}^\star$ is integrally closed in ${\mathbb Q}^\star$.
\item[(c)]
$({\mathbb Q}^{\star\star})^2\cap {\mathbb Q}^\star = ({\mathbb Q}^\star )^2$,
i.e. ${\mathbb Q}^\star$ is quadratically closed in ${\mathbb Q}^{\star\star}$.
\item[(d)]
If ${\mathbb Z}$ is diophantine in ${\mathbb Q}$
then ${\mathbb Z}^{\star\star}\cap {\mathbb Q}^\star = {\mathbb Z}^\star$
and ${\mathbb Q}^\star$ is algebraically closed in ${\mathbb Q}^{\star\star}$.
\item[(e)]
${\mathbb Q}$ is not model complete, i.e.,
there are $\mathbb{Q}^\star$ and $\mathbb{Q}^{\star\star}$ such that
${\mathbb Q}^\star$ is not existentially closed in ${\mathbb Q}^{\star\star}$.
\end{enumerate}
\end{proposition}
{\em Proof:}
(a) is an immediate consequence of our universal definition of ${\mathbb Z}$ in ${\mathbb Q}$.
The very same definition holds for ${\mathbb Z}^\star$ in ${\mathbb Q}^\star$
and for ${\mathbb Z}^{\star\star}$ in ${\mathbb Q}^{\star\star}$
(it is part of Th$({\mathbb Q})$ that all definitions of ${\mathbb Z}$ in ${\mathbb Q}$ are equivalent).
So if this universal formula holds for $x\in{\mathbb Z}^{\star\star}\cap{\mathbb Q}^\star$ in ${\mathbb Q}^{\star\star}$ it also holds in ${\mathbb Q}^\star$, i.e. $x\in {\mathbb Z}^\star$.

(b) is true because ${\mathbb Z}^{\star\star}$ is integrally closed in ${\mathbb Q}^{\star\star}$.

(c) follows from the fact that both being a square and, by Proposition \ref{new}(b), not being a square are diophantine in ${\mathbb Q}$.

(d) If ${\mathbb Z}$ is diophantine in ${\mathbb Q}$
then ${\mathbb Z}^{\star\star}\cap {\mathbb Q}^\star\supseteq {\mathbb Z}^\star$ and hence equality holds, by (a).

To show that then also ${\mathbb Q}^\star$ is algebraically closed in ${\mathbb Q}^{\star\star}$, let us observe that, for each $n\in {\mathbb N}$,
$$A_n:=\{ (a_0,\ldots ,a_{n-1})\in{\mathbb Z}^n\mid\exists x\in {\mathbb Z}\mbox{ with }x^n+a_{n-1}x^{n-1} + \ldots + a_0=0\}$$
is decidable: zeros of polynomials in one variable are bounded in terms of their coefficients, so one only has to check finitely many $x\in {\mathbb Z}$.
In particular, by (for short) Matiyasevich's Theorem,
there is an $\exists$-formula $\phi (t_0,\ldots ,t_{n-1})$ such that
$${\mathbb Z}\models\forall t_0\ldots t_{n-1}\left(\{\forall x
[ x^n+t_{n-1}x^{n-1}+\ldots +t_0\neq 0]\}\leftrightarrow \phi (t_0,\ldots ,t_{n-1})\right).$$
Since both $A_n$ and its complement in ${\mathbb Z}^n$ are diophantine in ${\mathbb Z}$, the same holds in ${\mathbb Q}$, by our assumption of ${\mathbb Z}$ being diophantine in ${\mathbb Q}$,
i.e. $A_n^{\star\star}\cap ({\mathbb Q}^\star)^n = A_n^\star$.
As any finite extension of ${\mathbb Q}^\star$ is generated by an integral primitive element this implies that ${\mathbb Q}^\star$ is relatively algebraically closed in ${\mathbb Q}^{\star\star}$.

(e) Choose a recursivley enumerable subset $A\subseteq {\mathbb Z}$ which is not decidable.
Then $B:={\mathbb Z}\setminus A$ is definable in ${\mathbb Z}$, and hence in ${\mathbb Q}$.
If $B$ were diophantine in ${\mathbb Q}$ it would be recursively enumerable. But then $A$ would be decidable: contradiction.

So not every definable subset of ${\mathbb Q}$ is diophantine in ${\mathbb Q}$,
and hence ${\mathbb Q}$ is not model complete.
Or, in other words, there are models ${\mathbb Q}^\star, {\mathbb Q}^{\star\star}$ of Th$({\mathbb Q})$ with ${\mathbb Q}^\star\subseteq {\mathbb Q}^{\star\star}$
where ${\mathbb Q}^\star$ is not existentially closed in ${\mathbb Q}^{\star\star}$.\qed
\medskip

We are confident that with similar methods as used in this paper one can show for an arbitrary prime $p$ that the unary predicate `$x\not\in {\mathbb Q}^p$' is also diophantine. This would imply that, in the setting of the Proposition, ${\mathbb Q}^\star$ is always radically closed in ${\mathbb Q}^{\star\star}$.
However, we have no bias towards an answer (let alone an answer) to the following (unconditional)
\begin{question}
For ${\mathbb Q}^\star\equiv {\mathbb Q}^{\star\star}\equiv {\mathbb Q}$ with
${\mathbb Q}^\star\subseteq {\mathbb Q}^{\star\star}$,
is ${\mathbb Q}^\star$ always algebraically closed in ${\mathbb Q}^{\star\star}$?
\end{question}

Mathematical Institute, 24-29 St Giles', Oxford OX1 3LB, UK\\
{\tt koenigsmann@maths.ox.ac.uk}
\end{document}